\documentclass{article}

\usepackage{amsmath}
\usepackage[latin1]{inputenc}
\usepackage[T1]{fontenc}
\usepackage[english]{babel}
\usepackage{enumerate}
\usepackage{amssymb}
\usepackage{amsthm}
\usepackage{cleveref}
\usepackage{graphicx}
\usepackage{url}

\newtheorem{Lemma}{Lemma}[section]
\newtheorem{Theorem}[Lemma]{Theorem}
\newtheorem*{Theorem*}{Theorem}
\newtheorem{Corollary}[Lemma]{Corollary}
\newtheorem*{Corollary*}{Corollary}
\newtheorem{Remark}[Lemma]{Remark}
\newtheorem{Example}[Lemma]{Example}

\title{The Hyperrigidity Conjecture for compact convex sets in $\mathbb{R}^2$}
\author{Marcel Scherer}
\date{}

\begin{document}

\maketitle

\begin{abstract}
We prove that for every compact, convex subset $K\subset\mathbb{R}^2$ the operator system $A(K)$, consisting of all continuous affine functions on $K$, is hyperrigid in the C*-algebra $C(\textup{ex}(K))$. In particular, this result implies that the weak and strong operator topologies coincide on the set
\begin{equation*}
    \{T\in\mathcal{B}(H);\ T\ \textup{normal}\ \textup{and}\ \sigma(T)\subset\textup{ex}(K)\}.
  \end{equation*}
Our approach relies on geometric properties of $K$ and generalizes previous results by Brown.
\end{abstract}

\section{Introduction}
\let\thefootnote\relax\footnotetext{ \date\\
2020 \textit{Mathematics Subject Classification.}Primary 46L05, 46L07, 47B15.\\
\ \textit{Key words and phrases.} operator system, hyperrigidity conjecture, unital completely positive, commutative $C^*$-algebra, convex set, affine function, polar parameterization.\\
 The author was partially supported by the Emmy Noether Program of the German Research Foundation (DFG Grant 466012782).}

An operator system $S$ that generates a C*-algebra $C^*(S)$ is called \textit{hyperrigid} if for every unital *-homomorphism $\pi:C^*(S)\to\mathcal{B}(H)$ has the \textit{unique extension property}, that is, the only unital completely positive map $\phi:C^*(S)\to\mathcal{B}(H)$ satisfying $\pi|_S=\phi|_S$ is $\pi$ itself. In 2011, Arveson conjectured in \cite{WA} that $S$ is hyperrigid if and only if the restrictions of all irreducible representations have the unique extension property. This problem continues to attract significant interest, see \cite{RC}, \cite{CT}, \cite{CS} and \cite{BB}. Although the conjecture has recently been disproven in \cite{BD}, its validity for operator systems in commutative C*-algebras remains unsolved, for partial results see \cite{DK}, \cite{PS} and \cite{WA}.\\
Among the contributions to the hyperrigidity conjecture, one noteworthy yet less well-known paper is Lawrence Brown's \cite{LB}, which shows that the operator systems $\langle 1,t,f\rangle$ are hyperrigid in $C([0,1])$ for strictly convex/concave $f$. There are already some results inspired by this, \cite{PS}, but we are going in a completely different direction. In this paper we show the following main \Cref{ThmMain}.
\begin{Theorem*}
If $K\subset\mathbb{R}^2$ is compact and convex, then $A(K)\subset C(\textup{ex}(K))$ is hyperrigid.
\label{ThmMain}
\end{Theorem*}
This generalizes the work of Brown since for strictly convex $f$ the operator system
  \begin{equation*}
    A(\textup{conv}(\{(t,f(t));\ t\in[0,1]\}))\subset C(\{(t,f(t));\ t\in[0,1]\})
  \end{equation*}
is completely order isomorphic to $\langle1,t,f\rangle$. The basic idea is similar to Brown's, but there are crucial differences. We cannot directly use left and right derivatives of functions in $A(K)$, instead we show in \Cref{Section2} that the polar parametrization of $K$ has left and right derivatives and use this to construct certain hyperplanes. Afterwards, we examine the function measuring the angles between these hyperplanes and collect several lemmas needed to prove the main theorem. At the beginning of \Cref{Section3}, we explain our approach in detail. Roughly speaking, we use a composition of a translation and rotation to separate connected disjoint subsets of $\partial K$. This leads us to two values which are given by basic geometric properties of $K$ and are studied in \Cref{Section3.1} respectively \Cref{Section3.2}. Afterwards we are finally ready to prove \Cref{ThmMain} and the following corollary in \Cref{Section3.3}.
  \begin{Corollary*}
For every compact, convex set $K\subset\mathbb{R}^2$, the weak and strong operator topologies coincide on the set
  \begin{equation*}
    \{T\in\mathcal{B}(H);\ T\ \textup{normal}\ \textup{and}\ \sigma(T)\subset\textup{ex}(K)\}.
  \end{equation*}
\end{Corollary*}
It is noteworthy that parts of our approach even works for compact convex sets in $\mathbb{R}^d$ for $d>2$, however, the main problem seems to be in \Cref{EqI} since a arbitrary path in $\partial K$ does in general not have finite length.

\section*{Acknowledgement}
I sincerely thank Michael Hartz, Malte Gerhold, Matthew Kennedy and Orr Shalit for listening to my ideas and valuable contributions to enhancing the clarity of the proofs presented in the paper.

\section{Preliminaries}
\label{Section2}
The following lemma is standard, and parts of it can be found in \cite{DD} or other classical books on convex analysis.
\begin{Lemma}
Let $I\subset\mathbb{R}$ be an open interval, and let $f:I\to\mathbb{R}$ be a convex function. Then
  \begin{enumerate}[(i)]
    \item the right $f'_+(t)$ and left $f'_-(t)$ derivatives exist for every $t\in I$,
    \item $f'_+$ is right continuous, $f'_-$ is left continuous,
    \item $\lim\limits_{s\to t, s<t\in I} f'_+(s)=f'_-(t)$ and $\lim\limits_{s\to t, s>t\in I} f'_-(s)=f'_+(t)$.
  \end{enumerate}
\label{LemBas}
\end{Lemma}

\underline{Proof}:\\
Part $(i)$ is well known and can be proven by using \cite[Lemma 16.5.8]{DD}, the Secant Lemma:
  \begin{equation}
    \frac{f(x)-f(a)}{x-a}\le\frac{f(b)-f(a)}{b-a}\le\frac{f(b)-f(x)}{b-x}
\label{Eq1}
  \end{equation}
for all $a<x<b\in I$.\\
For part $(ii)$, let $\epsilon>0$, $(s_n)_n$ in $I$ such that $s_n\to t \in I, t<s_n$ and $\tilde s\in I$ such that
  \begin{equation*}
    0\le\frac{f(s)-f(t)}{s-t}-f'_+(t)<\epsilon
  \end{equation*}
for every $t<s\le\tilde s$. From \Cref{Eq1}, we conclude that
  \begin{equation*}
      0\le \limsup_{n\to\infty}f'_+(s_n)-f'_+(t)\le\limsup_{n\to\infty}\frac{f(s_n)-f(\tilde s)}{s_n-\tilde s}-f'_+(t)=\frac{f(t)-f(\tilde s)}{t-\tilde s}-f'_+(t)<\epsilon,
  \end{equation*}
where we also used that a convex function is automatically continuous \cite[Theorem 16.5.9]{DD}. Therefore  $f'_+$ is right continuous and similarly we obtain that $f'_-$ is left continuous.\\
For the last part, let $s_n<t$ such that $s_n\to t$. The Secant Lemma shows that
  \begin{equation*}
    f'_-(a)\le f'_+(a)\le f'_-(b)
  \end{equation*}
for $a<b\in I$ (\cite[Theorem 16.5.9]{DD}). Hence the lower continuity of $f'_-$ implies that
  \begin{equation*}
    \lim_{n\to\infty}f'_-(s_n)=\lim_{n\to\infty}f'_+(s_n)=f'_-(t).
  \end{equation*}
The second part of $(iii)$ follows similarly. $\hfill\square$\\
\\
Let $K$ be a non-empty compact convex subset in $\mathbb{R}^2$ with $0\in \textup{Int}(K)$. The \textit{Minkowski functional}, given by
  \begin{equation*}
    f_K:\mathbb{R}^2\to\mathbb{R}, x\mapsto\inf\{r\in\mathbb{R}: r>0\textup{\ and\ }x\in rK\},
  \end{equation*}
is known to be a convex, continuous and nonnegative homogeneous function. Using this functional, the \textit{polar parameterization} of the boundary of $K$ is defined by
  \begin{equation*}
    p:\mathbb{R}\to \partial K, t\mapsto \frac{(\cos(t),\sin(t))}{f_K((\cos(t),\sin(t)))}.
  \end{equation*}
It is obvious that this is a $2\pi$-periodic function and in the next lemma we show that $p$ has the similar properties to those we saw in \Cref{LemBas}.\\
Throughout the paper, we identify $\mathbb{R}^2$ and $\mathbb{C}$ because this makes the proofs far more readable. However, it should be noted that the canonical scalar product we will be working with is only $\mathbb{R}$-linear. 

\begin{Lemma}
The following assertions hold:
 \begin{enumerate}
    \item The left and right derivatives $p'_-(\cdot)$ and $p'_+(\cdot)$ exist.
    \item $p'_-$ is left continuous and $p'_+$ is right continuous.
    \item $\lim\limits_{s\to t, s<t}p'_+(s)=p'_-(t)$ and $\lim\limits_{s\to t, s>t}p'_-(s)=p'_+(t)$.
  \end{enumerate}
  \label{LemCont}
\end{Lemma}

\underline{Proof}:\\
First, we show the lemma for $p$, restricted to the interval  $I=(-\pi/2,\pi/2)$. Let $h(t)=1/\cos(t)$. The function
  \begin{equation*}
    \tilde p(t)=f_K(h(t)e^{it})=f_K(1+i\tan(t))=f_K(1+i(\cdot))\circ \tan(t)
  \end{equation*}
is the composition of the convex function $f_K(1+i(\cdot))$ and $\tan(\cdot)$. By \Cref{LemBas}, we obtain that $f_K(1+i(\cdot))$ has properties $(i), (ii)$ and $(iii)$ from \Cref{LemBas}, and since $\tan(\cdot)$ is monotonically increasing and continuously differentiable, also the function $\tilde p$ has properties $(i), (ii), (iii)$ from \Cref{LemBas}. The nonnegative homogeneity of $f_K$ implies that
  \begin{equation*}
    \begin{split}
      &h(a)\frac{f_K(e^{ia})-f_K(e^{is})}{a-s}-\frac{\tilde p(a)-\tilde p(s)}{a-s}\\
       &=\frac{h(a)(f_K(e^{ia})-f_K(e^{is}))-f_K(h(a)e^{ia})+f_K(h(s)e^{is})}{a-s}\\
       &=f_K(e^{is})\frac{h(s)-h(a)}{a-s}
    \end{split}
  \end{equation*}
for all $a\in I$, so $f_K(e^{i\cdot})$ is right and left differentiable in $a$ with
  \begin{equation*}
    \begin{split}
       f'_{K+}(e^{ia})=\frac{1}{h(a)}(\tilde p'_+(a)-f_K(e^{ia})h'(a))\\
       f'_{K-}(e^{ia})=\frac{1}{h(a)}(\tilde p'_-(a)-f_K(e^{ia})h'(a)).
    \end{split}
  \end{equation*}
Therefore $f_K$ has properties $(ii), (iii)$ from \Cref{LemBas} and we can immediately conclude that $p$ restricted to $(-\pi/2,\pi/2)$ has property $1., 2.$ and $3.$.\\
To finish the proof, one does the same as above with $I=(0,\pi)$ and $h(t)=1/\sin(t)$, $I=(-\pi,0)$ and $h(t)=-1/\sin(t)$, $I=(\pi/2,3\pi/2)$ and $h(t)=-1/\cos(t)$. $\hfill\square$\\
\\

\begin{figure}[h]
    \centering
    \includegraphics[width=0.6\textwidth]{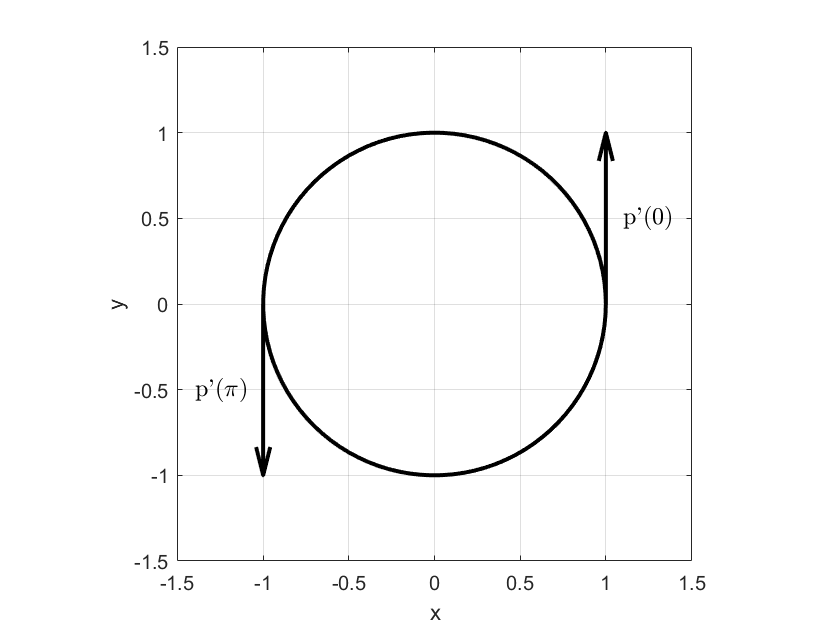} 
    \caption{$\partial\mathbb{D}$ with $p'(0)$ and $p'(\pi)$}
    \label{Fig1}
\end{figure}

\begin{Example}
Let $K$ be the unit disc. Then the polar parametrization is given by
  \begin{equation*}
    p(t)=(\cos(t),\sin(t))
  \end{equation*}
The derivative is of course $p'(t)=(-\sin(t),\cos(t))$. \Cref{Fig1} visualizes $\partial K$ and $p'(0)$ and $p'(\pi)$. It should be noted at this point that the direction of the vectors $p'(t)$ will be important later on.
\end{Example}

For $t\in\mathbb{R}$, we define the subderivatives $\partial^-_{p(t)}=\{p(t)+sp'_-(t), s\in\mathbb{R}\}$ and $\partial^+_{p(t)}=\{p(t)+sp'_+(t), s\in\mathbb{R}\}$. In the next theorem, we will show that $\partial^+_{p(t)}$ and $\partial^-_{p(t)}$ are supporting hyperplanes. Recall that a supporting hyperplane of $K$ is a hyperplane such that $K$ is entirely contained in one of the two closed half-spaces bounded by the hyperplane and $K$ intersects the hyperplane.

\begin{figure}[h]
    \centering
    \includegraphics[width=0.6\textwidth]{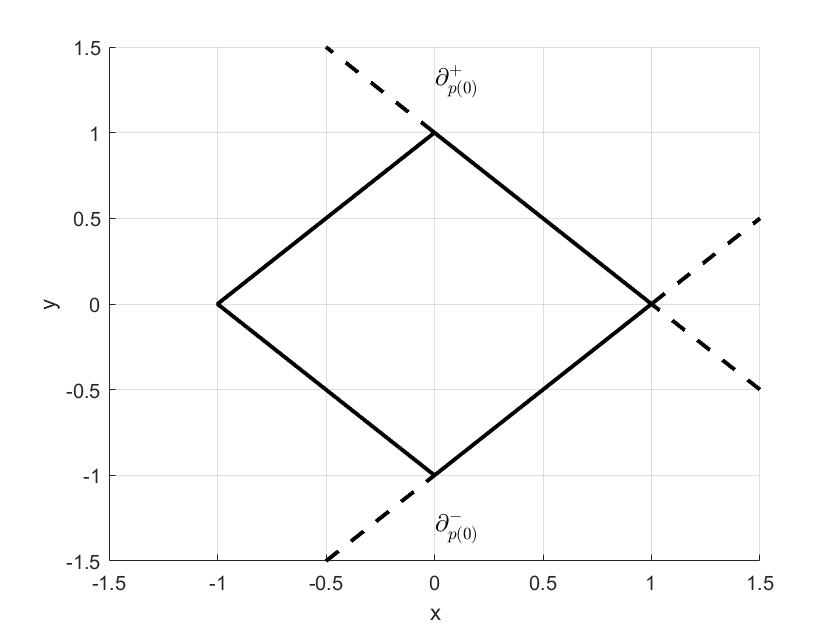} 
    \caption{an example of $\partial K$ with $\partial^+_{p(0)}$ and $\partial^-_{p(0)}$}
    \label{Fig2}
\end{figure}

\begin{Example}
Let $K$ be the convex hull of $\{1, -1, i, -i\}$. Then, the polar parametrization $p$ is given by
  \begin{equation*}
    p(t)=(0,1)+\frac{1}{1+\tan(t)}(1,-1)
  \end{equation*}
on $(0,\pi/2)$ and 
  \begin{equation*}
    p(t)=(1,0)+\frac{1}{1-\tan(t)}(1,1)
  \end{equation*}
on $(3\pi/2,2\pi)$. Therefore $p'_+(0)=(-1,1)$, $p'_-(0)=(1,1)$ and 
  \begin{equation*}
    \begin{split}
      &\partial^+_{p(0)}=\{(0,1)+s(-1,1);\ s\in\mathbb{R}\},\\
      &\partial^-_{p(0)}=\{(0,1)+s(1,1);\ s\in\mathbb{R}\}.
    \end{split}
  \end{equation*}
For a better visualization, this example is illustrated in \Cref{Fig2}.
\end{Example}

\begin{Theorem}
For every $x\in\partial K$, $\partial^-_{x}$ and $\partial^+_x$ are supporting hyperplanes of $K$.
\label{ThmHyp}
\end{Theorem}

\underline{Proof}:\\
The first thing to check is that $\partial^-_{x}$ and $\partial^+_x$ are actually hyperplanes and not just points. For this purpose, it suffices to verify that $p'_-(t)\neq0$ and $p'_+(t)\neq0$. Note that we saw in the proof of \Cref{LemCont} that $f_K(e^{i\cdot})$ is left and right differentiable. Thus, the required inequalities can be directly deduced from the following inequalities
  \begin{equation*}
    \begin{split}
      \lim_{s\to t, s<t}\left|\frac{p(t)-p(s)}{t-s}\right|&=\lim_{s\to t, s<t}\frac{1}{f_K(e^{is})f_K(e^{it})}\left|\frac{f_K(e^{is})e^{it}-f_K(e^{it})e^{is}}{t-s}\right|\\
      &=\frac{1}{f_K(e^{it})^2}\lim_{s\to t, s<t}\left|\frac{f_K(e^{is})e^{it}-f_K(e^{is})e^{is}+f_K(e^{is})e^{is}-f_K(e^{it})e^{is})}{t-s}\right|\\
      &=\frac{1}{f_K(e^{it})^2}\lim_{s\to t, s<t}\left|f_K(e^{is})\frac{e^{it}-e^{is}}{t-s}+e^{is}\frac{f_K(e^{is})-f_K(e^{it})}{t-s}\right|\\
      &=\frac{1}{f_K(e^{it})^2}\left|ie^{it}f_K(e^{it})+e^{it}f_{K}(e^{i\cdot})'_-(t)\right|\\
      &=\frac{1}{f_K(e^{it})^2}\left|if_K(e^{it})+f_{K}(e^{i\cdot})'_-(t)\right|\\
      &\ge\frac{1}{f_K(e^{it})}>0
    \end{split}
  \end{equation*}
and $|p'_+(t)|\ge\frac{1}{f_K(e^{it})^2}|if_K(e^{it})+f_{K}(e^{i\cdot})'_+(t)|\ge\frac{1}{f_K(e^{it})}>0$.\\
\\
Now that this has been shown, it is time to prove that $\partial^-_x$ and $\partial^+_x$ are supporting hyperplanes. Let $x=p(t)\in\partial K$ with $t\in [0,2\pi)$ and $(t_n)_n$ a sequence in $(0,2\pi)$ such that $t_n\to t$ and $t_n>t$. The line through the points $p(t)$ and $p(t_n)$ splits $\mathbb{R}^2$ into two closed half planes
  \begin{equation*}
    \begin{split}
      & H^n_1=\left\{z\in\mathbb{R}^2;\ \left\langle z-p(t),\frac{i(p(t)-p(t_n))}{t-t_n}\right\rangle\ge0\right\},\\ 
      &H^n_2=\left\{z\in\mathbb{R}^2;\ \left\langle z-p(t),\frac{i(p(t)-p(t_n))}{t-t_n}\right\rangle\le0\right\}.
    \end{split}
  \end{equation*}
 The continuity of the polar parameterization implies that
  \begin{equation*}
    \{p(s);\ s\in[0,2\pi)\setminus(t,t_n)\}
  \end{equation*}
is either contained in $H^n_1$ or $H^n_2$, and by the Pigeonhole principle, we have either for $i=1$ or $i=2$ that 
  \begin{equation*}
    p([0,2\pi)\setminus(t,t-n))\subset H^n_i
  \end{equation*}
for infinitely many $n$. On the other hand the two closed half spaces bounded by $\partial^+_x$ are given by 
  \begin{equation*}
    \begin{split}
      &\tilde H_1=\{z\in\mathbb{R}^2;\ \langle z-p(t),ip'_+(t)\rangle\ge0\},\\ 
      &\tilde H_2=\{z\in\mathbb{R}^2;\ \langle z-p(t),ip'_+(t)\rangle\le0\}.
    \end{split}
  \end{equation*}
Thus, since $p_+'(t)=\lim_{n\to\infty}\frac{p(t)-p(t_n)}{t-t_n}$, either $\tilde H_1$ or $\tilde H_2$ contains
\begin{equation*}
    \{p(s);\ s\in[0,2\pi)\}.
  \end{equation*} 
Therefore $\partial^+_x$ is a supporting hyperplane and analogously one shows that $\partial^-_x$ is also a supporting hyperplane. $\hfill\square$\\
\\

\begin{Remark}
If there is only one supporting hyperplane $\partial_x$ for a point $x\in\partial K$, then \Cref{ThmHyp} yields
  \begin{equation*}
    \partial^+_x=\partial^-_x=\partial_x.
  \end{equation*}
In particular, let $[a,b]=I\subset[0,2\pi)$ such that $\{p(t);\ t\in[a,b]\}$ is contained in a face of $K$. For every $s\in(a,b)$, there is only one supporting hyperplane through $p(s)$ and it is given by the line through $p(a)$ and $p(b)$, denoted by $[p(a):p(b)]$. Thus 
  \begin{equation*}
    \partial^+_{p(s)}=\partial^-_{p(s)}=[p(a):p(b)]
  \end{equation*}
for all $s\in(a,b)$ and by \Cref{LemCont} we also obtain that
  \begin{equation*}
    \partial^+_{p(a)}=[p(a):p(b)]=\partial^-_{p(b)}.
  \end{equation*}
\label{RemHyp}
\end{Remark}

The next step is to define a function $\sphericalangle(\cdot,\cdot):[0,2\pi]^2\to[0,\pi]$ that measures the angle between two subderivatives. To do this, we need the following two pieces of notation. For three different points $x, y, z\in\mathbb{R}^2$, denote by $\triangle(x,y,z)$ the triangle with corners $x, y, z$ and by $\angle(x, y, z)$ the angle in the point $y$. Now define for $s\le t$
  \begin{equation*}
    \sphericalangle(s,t)=\begin{cases}
         0 & \textup{if}\ \partial^+_{p(s)}\cap\partial^-_{p(t)}=\emptyset\\
        \pi & \textup{if}\ \partial^+_{p(s)}=\partial^-_{p(t)}\\
        \cos^{-1}\left(\frac{\langle p'_-(t),-p'_+(s)\rangle}{|p'_-(t)||p'_+(s)|}\right) & \textup{if}\ \partial^+_{p(s)}\cap\partial^-_{p(t)}=\{z\} \ \textup{and}\\ & \{p(l);\ l\in[s,t]\}\subset\triangle(p(s),z,p(t))\\
       0 & \textup{else}
      \end{cases}
  \end{equation*}
and for $s>t$, $\sphericalangle(s,t)=\sphericalangle(t,s)$. One should note that in the third case  
  \begin{equation}
   \sphericalangle(s,t)=\angle(p(s), z, p(t))
  \label{Eq2}
  \end{equation}
 for $s\neq t$.
\\

\begin{Remark}
A careful reader should ask themselves why $\angle(p(s), z, p(t))$ is not given by
  \begin{equation*}
    \cos^{-1}\left(\frac{\langle p'_-(t),p'_+(s)\rangle}{|p'_-(t)||p'_+(s)|}\right).
  \end{equation*}
To prove \Cref{Eq2}, it suffices to show that for every $t\in[0,2\pi)$ the set $K$ is to the left of $\{p(t)+sp'_-(t);\ s\in\mathbb{R}\}$ resp. $\{p(t)+sp'_+(t);\ s\in\mathbb{R}\}$. Since $0\in\textup{Int}(K)$, this directly follows from
  \begin{equation*}
    \begin{split}
      &\langle 0-p(t),ip'_+(t)\rangle=\lim_{n\to\infty, t_n>t}\left\langle-p(t),i\frac{p(t)-p(t_n)}{t-t_n}\right\rangle\\
      &=\lim_{n\to\infty, t_n>t}\frac{1}{f_K(e^{it})f_K(e^{it_n})(t-t_n)}\langle (\cos(t),\sin(t)),(-\sin(t_n),\cos(t_n))\rangle\\
      &=\lim_{n\to\infty, t_n>t}\frac{-1}{f_K(e^{it})f_K(e^{it_n})(t_n-t)}(-\cos(t)\sin(t_n)+\sin(t)\cos(t_n))\ge0,
    \end{split}    
  \end{equation*}
whereby the last inequality follows for $t\neq \pi/2$ and $t\neq3/2\pi$ and $t_n$ sufficiently close to $t$ from 
  \begin{equation*}
    \frac{\sin(t)}{\cos(t)}=\tan(t)\le\tan(t_n)=\frac{\sin(t_n)}{\cos(t_n)},
  \end{equation*}
and $\sin(t)\cos(t_n)\le0$ for $t=\pi/2$ or $t=3/2\pi$.\\
Similarly one obtains that $\langle 0-p(t),ip'_-(t)\rangle\ge0$. 
\label{RemLeft}
  \end{Remark}

\begin{Lemma}
If $s_n<t_n$ in $[0,2\pi]$ such that $s_n, t_n\to t$ and $s_n< t< t_n$, then
  \begin{enumerate}[(a)]
    \item $\lim_{n\to\infty}\sphericalangle(s_n,t_n)=\sphericalangle(t,t)$.
  \end{enumerate}
If $t\le s_n<t_n$ and $s_n, t_n\to t$, then
  \begin{enumerate}[(b)]
    \item $\lim_{n\to\infty}\sphericalangle(s_n,t_n)=\pi$,
  \end{enumerate}
and if $t_n<s_n\le t$ and $s_n, t_n\to t$, then
  \begin{enumerate}[(c)]
    \item $\lim_{n\to\infty}\sphericalangle(t_n,s_n)=\pi$,
   \end{enumerate}
\label{LemCon}
\end{Lemma}

\underline{Proof}:\\
First we show that for every $x\in\mathbb{R}$ there is a $\epsilon>0$ such that for all $s<t\in B_\epsilon(x)$ we have that either $\partial^+_{p(s)}=\partial^-_{p(t)}$ or $\partial^+_{p(s)}\cap\partial^-_{p(t)}=\{z\}$ and $\{p(l);\ l\in[s,t]\}\subset\triangle(p(s),z,p(t))$. If we assume the opposite, then for every $n\in\mathbb{N}$ there are $s_n<t_n\in B_{1/n}(x)$ such that $\partial^+_{p(s_n)}\cap\partial^-_{p(t_n)}=\emptyset$ or $\partial^+_{p(s_n)}\cap\partial^-_{p(t_n)}=\{z_n\}$ and $\partial K\setminus\{p(l);\ l\in[s_n,t_n]\}\subset\triangle(p(s_n),z_n,p(t_n))$. Note that in the first case, $K$ is between the two parallel lines $\partial^+_{p(s_n)}$ and $\partial^-_{p(t_n)}$ and in the second case, 
  \begin{equation*}
    \{p(t);\ t\in[x-\pi,x+\pi]\setminus(s_n,t_n)\}\subset\triangle(p(s_n),z_n,p(t_n)).
  \end{equation*}
Because $s_n, t_n\to x$ for $n\to\infty$ and because $p$ is continuous, we obtain the contradiction that $0\in\textup{Int}(K)$. \\
The lemma now follows from \Cref{Eq2} and \Cref{LemCont}. $\hfill\square$

\begin{Lemma}
Let $a<b\in[0,2\pi)$ and $\epsilon>0$, then there are $a=t_0<t_1<\dots<t_m=b$ such that 
  \begin{equation*}
    \sphericalangle(t_n,t_{n+1})\ge\pi-\epsilon.
  \end{equation*}
\label{LemPar}
\end{Lemma}

\underline{Proof}:\\
First, note that the set 
  \begin{equation*}
    \{x\in(a,b);\ \sphericalangle(x,x)\le\pi-\epsilon\}
  \end{equation*}
contains only finitely many points $s_1<s_2<\dots<s_k$. If we assume the opposite, then for every $n$ there is a convex polygon such that each angle in each corner is smaller than $\pi-\epsilon$. However, the sum of the angles is also given by $(n-2)\pi$, which is a contradiction for $n$ big enough.\\ 
Let $s_0=a, s_{k+1}=b$. Next we show that there is a $\delta>0$ such that 
  \begin{equation*}
    \sphericalangle(x,y)\ge\pi-\epsilon
  \end{equation*}
for all $x, y\in[s_n,s_{n+1}], n=0,\dots, k$ and $|x-y|<\delta$. Assuming the opposite yields a $m\in\{0,\dots,k\}$ and $x_n<y_n\in[s_m,s_{m+1}]$ such that $|x_n-y_n|<1/n$ and $\sphericalangle(x_n,y_n)\le\pi-\epsilon$. Without loss of generality let $(x_n)_n, (y_n)_n$ converge to a $z\in[s_m,s_{m+1}]$. By \Cref{LemCon} it holds that either 
  \begin{equation*}
    \lim_{n\to\infty}\sphericalangle(x_n,y_n)=\pi
  \end{equation*}
if $z\le x_n$ or $y_n\le z$, or 
 \begin{equation*}
    \lim_{n\to\infty}\sphericalangle(x_n,y_n)=\sphericalangle(z,z)
  \end{equation*}
if $x_n<z<y_n$. In the latter case, it holds that $\sphericalangle(z,z)>\pi-\epsilon$ by choice of $s_1,\dots, s_k$, contradicting $\sphericalangle(x_n,y_n)\le \pi-\epsilon$.\\
Thus we obtain the desired partition as a refinement of $s_0<\dots<s_{k+1}$. $\hfill\square$

\section{Hyperrigidity of $A(K)$}
\label{Section3}
Let $K\subset \mathbb{R}^2$ be compact and convex. If $\textup{Int}(K)=\emptyset$, then $K$ is merely an interval and $A(K)=C(\textup{ex}(K))$. Therefore, $A(K)$ is clearly hyperrigid in $C(\textup{ex}(K))$. If $x\in\textup{Int}(K)$, then define $\tilde K=\{y-x;\ y\in K\}$ and 
  \begin{equation*}
    j:C(\textup{ex}(K))\to C(\textup{ex}(\tilde K)), f\mapsto f(\cdot+x),
  \end{equation*}
where we do not have to close $\textup{ex}(K)$ since the extreme points of closed convex sets in $\mathbb{R}^2$ are always closed. It is clear that $j$ is a unital *-isomorphism that maps $A(K)$ onto $A(\tilde K)$. Therefore showing that $A(K)$ is hyperrigid in $C(\textup{ex}(K))$ is equivalent to showing that $A(\tilde K)$ is hyperrigid in $C(\textup{ex}(\tilde K))$. Note that $0\in\textup{Int}(\tilde K)$ by construction and therefore we will always assume that $0\in\textup{Int}(K)$. This ensures that we can use the polar parametrization $p$. \\
Let $H$ be a Hilbert space, $\pi: C(\textup{ex}(K))\to\mathcal{B}(H)$ be a unital *-homomorphism and $\phi:C(\textup{ex}(K))\to\mathcal{B}(H)$ be a u.c.p. map such that $\pi=\phi$ on $A(K)$. For $F\subset\mathbb{R}^2$ Borel, we define
  \begin{equation*}
    \begin{split}
      &\phi(\chi_{F})=\mu(F\cap\textup{ex}(K)),\\ 
      &\pi(\chi_{F})=\nu(F\cap\textup{ex}(K)), 
    \end{split}
  \end{equation*}
where $\mu$ respectively $\nu$ is the positive operator valued measure corresponding to $\phi$ respectively $\pi$.\\
We now briefly outline the general strategy behind the proof of the main theorem. Like in Brown's proof, we want to show that $\pi(\chi_{p(I)})\phi(\chi_{p(J)})=0$ for all disjoint closed intervals $I, J\subset[0,2\pi)$. Form there, it won't be too far to see that $\pi=\phi$. To achieve $\pi(\chi_{p(I)})\phi(\chi_{p(J)})=0$, we will cover $I$ with smaller disjoint intervals $I_i$ so that we can use \cite[Corollary 1.2]{LB}, which will play a crucial role in our proof. For the convenience of the reader, we will state this corollary in the following lemma. 

\begin{Lemma}
Let $H=H_1\oplus\dots, H_n$ be a direct sum of Hilbert spaces, $A$ a positive operator on $H$ and $P_{i}$ the projection onto $H_i$. Then
  \begin{equation*}
    \|A\|\le\sum_{i=1}^n\|P_iAP_i\|.
  \end{equation*}
\label{LemBrown}
\end{Lemma}

In our setting, we choose $A=\pi(\chi_{p(I)})\phi(\chi_{p(J)})\pi(\chi_{p(I)})$ and $H_i=\pi(\chi_{p(I_i)})(H)$, whereby $(I_i)_{i=1}^n$ is a family of disjoint Borel sets such that $I=\bigcup_{i=1}^nI_i$. Thus the lemma yields that
  \begin{equation*}
    \|\pi(\chi_{p(I)})\phi(\chi_{p(J)})\pi(\chi_{p(I)})\|\le\sum_{i=1}^n\|\pi(\chi_{p(I_i)})\phi(\chi_{p(J)})\pi(\chi_{p(I_i)})\|.
  \end{equation*}
Hence we only need a good upper estimate for $\|\pi(\chi_{p(I_i)})\phi(\chi_{p(J)})\pi(\chi_{p(I_i)})\|$. Let us drop the index $i$ for a moment and take $f\in A(K)$ with $\chi_{p(J)}\le f$. Then the positivity of $\phi$ yields
  \begin{equation}
    \pi(\chi_{p(I)})\phi(\chi_{p(J)})\pi(\chi_{p(I)})\le\pi(\chi_{p(I)})\phi(f)\pi(\chi_{p(I)})=\pi(f\chi_{p(I)}).
\label{EqPhi}
  \end{equation}
So the next step is to construct a suitable $f$. Let $I=[a,b]\subset[0,2\pi)$ with $a<b$ and $g_I$ be the unique composition of a translation and a rotation such that $\textup{Im}(g(p(a)))=\textup{Im}(g(p(b)))=0$ and $\textup{Re}(g(p(a)))<\textup{Re}(g(p(b)))=0$ and assume that $\textup{inf}_{p(J)}|\textup{Im}(g_I)|>0$. Then the function
  \begin{equation*}
    \tilde g_I=\frac{\textup{Im}(g_I)+\|\textup{Im}(g_I)\|_{p(I),\infty}}{\textup{inf}_{p(J)}|\textup{Im}(g_I)|}
  \end{equation*}
is an affine function which, when plugged into \Cref{EqPhi}, yields that
  \begin{equation*}
    \|\pi(\chi_{p(I)})\phi(\chi_{p(J)})\pi(\chi_{p(I)})\|\le\frac{\|\textup{Im}(g_I)\|_{p(I),\infty}}{\textup{inf}_{p(J)}|\textup{Im}(g_I)|},
  \end{equation*}
since $\pi(\textup{Im}(g_I)\chi_{p(I)})\le0$. Therefore we establish an upper bound for $\|\textup{Im}(g_I)\|_{p(I),\infty}$ in \Cref{Section3.1} and a lower bound for $\textup{inf}_{p(J)}|\textup{Im}(g_I)|$  in \Cref{Section3.2}.

\subsection{The upper bound for $\|\textup{Im}(g_I)\|_{p(I), \infty}$}
\label{Section3.1}
\begin{figure}[h]
    \centering
    \includegraphics[width=0.6\textwidth]{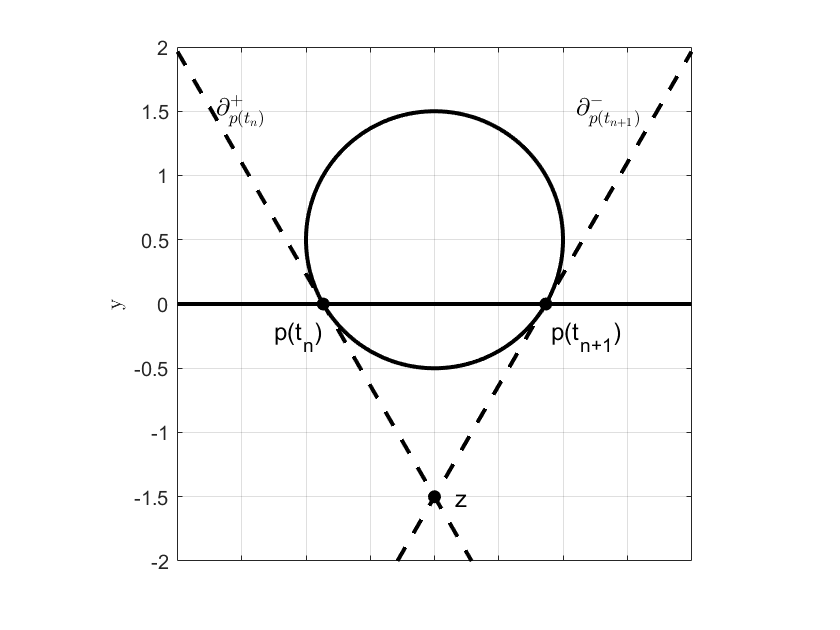} 
    \caption{an illustration for the proof of \Cref{LemUpp}}
    \label{Fig3}
\end{figure}

\begin{Lemma}
Let $I=[a,b], J\subset[0,2\pi)$ be closed disjoint intervals with $a<b$. Then, for every $\epsilon>0$, there exists a partition $a=t_0<t_1<\dots<t_m=b$ such that
  \begin{equation*}
    \|\textup{Im}(g_{[t_n,t_{n+1}]})\|_{p([t_n,t_{n+1}]),\infty}\le\epsilon|p(t_n)-p(t_{n+1})|
  \end{equation*}
for every $n=0,\dots,m-1$.
\label{LemUpp}
\end{Lemma}

The idea of the proof is illustrated in \Cref{Fig3}. \\
\underline{Proof}:\\
Let $\pi/2>\epsilon>0$. By \Cref{LemBas} there is a partition $a=t_0<t_1<\dots<t_m=b$ such that 
  \begin{equation*}
    \sphericalangle(t_n,t_{n+1})\ge\pi-\epsilon
  \end{equation*}
for every $0\le n\le m-1$. Fix $n$. For the sake of clarity, let $g=g_{[t_n,t_{n+1}]}$ and $I_n=[t_n,t_{n+1}]$. If $\partial^+_{p(t_n)}=\partial^-_{p(t_{n+1})}$, then $\|\textup{Im}(g)\|_{p(I_n), \infty}=0$. Therefore, we can assume that $\{z\}=\partial^+_{p(t_n)}\cap\partial^-_{p(t_{n+1})}$ for some $z\in\mathbb{R}^2$. By the definition of $\sphericalangle(\cdot,\cdot)$ and choice of the partition we have that
  \begin{equation*}
     p(I_n)\subset \triangle(p(t_n), z, p(t_{n+1}))
  \end{equation*}
and therefore $\|\textup{Im}(g)\|_{p(I_n),\infty}\le|\textup{Im}(g(z))|$. It also holds that
  \begin{equation*}
    \sphericalangle(t_n,t_{n+1})\ge\pi-\epsilon>\pi/2
  \end{equation*}
and thus $\angle(p(t_{n+1}),p(t_n),z)<\pi/2$ and $\angle(z,p(t_{n+1}),p(t_n))<\pi/2$. This yields that $\textup{Re}(g(p(t_{n+1})))\le\textup{Re}(g(z))\le\textup{Re}(g(p(t_n)))$ and everything together with basic geometry results in 
  \begin{equation*}
    \begin{split}
      \frac{\|\textup{Im}(g)\|_{p(I_n),\infty}}{|p(t_n)-p(t_{n+1})|}&=\frac{\|\textup{Im}(g)\|_{p(I_n),\infty}}{|\textup{Re}(g(p(t_n)))-\textup{Re}(g(p(t_{n+1})))|}\\
      &\le\frac{|\textup{Im}(g(z))|}{|\textup{Re}(g(p(t_n)))-\textup{Re}(g(z))|}\\
      &=\tan(\angle(g(p(t_{n+1})),g(p(t_n)),g(z))\\
      &=\tan(\angle(p(t_{n+1}),p(t_n),z)\\
      &\le\tan(\pi-\sphericalangle(t_n,t_{n+1}))\le \tan(\epsilon).
    \end{split}
  \end{equation*}
Thus we obtain the lemma by making $\epsilon$ sufficiently small. $\hfill\square$\\

\subsection{The lower bound for $\textup{inf}_{p(J)}|\textup{Im}(g_I)|$}
\label{Section3.2}

For two sets $A, B\subset \mathbb{R}^2$ we define the distance between $A$ and $B$ as
  \begin{equation*}
    \textup{dist}(A,B)=\textup{inf}\{|x-y|;\ x\in A, y\in B\}.
  \end{equation*}
\begin{Lemma}
Let $I, J$ be closed disjoint intervals in $[0,2\pi)$ and $I=[\alpha,\beta]$ with $\alpha<\beta$. Then
  \begin{equation*}
    \textup{inf}_{p(J)}|\textup{Im}(g_{I})|=\textup{dist}(g_I(J),\mathbb{R}\times\{0\})\ge\textup{dist}(p(J),\partial^+_{p(\alpha)}\cup\partial^-_{p(\beta)}\cup[p(\alpha):p(\beta)]),
\end{equation*}
where $[p(\alpha):p(\beta)]$ denotes the hyperplane through $p(\alpha)$ and $p(\beta)$. 
\label{LemDis1}
\end{Lemma}

\underline{Proof}:\\
First note that
  \begin{equation}
    \begin{split}
      \textup{inf}_{p(J)}|\textup{Im}(g_I)|&=\inf\{|x-g_I(y)|;\ x\in\mathbb{R}\times\{0\}, y\in p(J)\}\\
      &=\textup{dist}(g_I(p(J)),\mathbb{R}\times\{0\}).
\label{EqDist}
    \end{split}
  \end{equation}
Now observe that $g_{I}([p(\alpha):p(\beta)])=\mathbb{R}\times\{0\}$ since $g_I(p(\alpha)), g_I(p(\beta))\in\mathbb{R}\times\{0\}$ by construction of $g_I$, and that since we chose $g_I$ to be a composition of a translation and rotation, we have that $|g_I(x)-g_I(y)|=|x-y|$ for all $x, y\in\mathbb{R}^2$. These two observations combined with \Cref{EqDist} yield the lemma. $\hfill\square$\\
\\
\begin{figure}[h]
    \centering
    \includegraphics[width=0.6\textwidth]{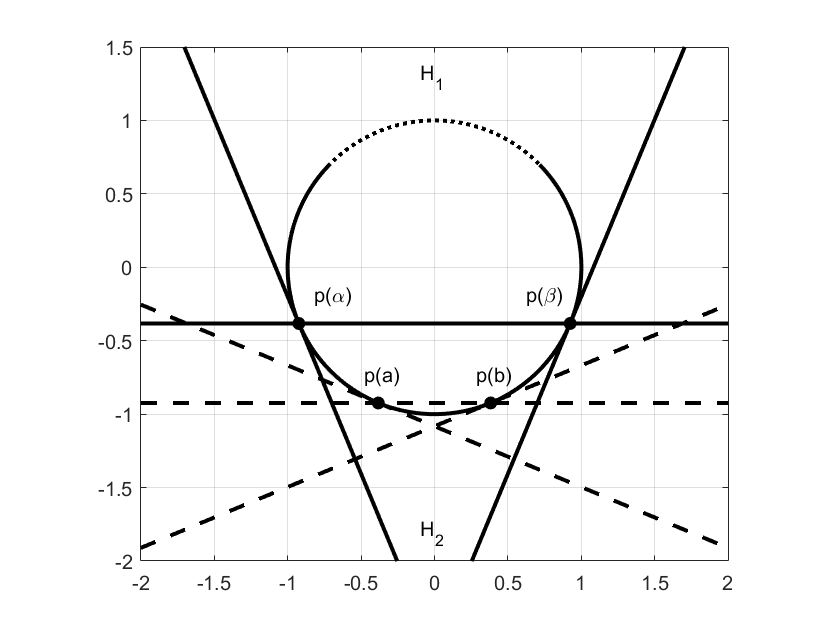} 
    \caption{$\partial\mathbb{D}$ with $\alpha=9/8\pi$, $a=11/8\pi$, $b=13/8\pi$ and $\beta=15/8\pi$. The hyperplanes $\partial^+_{p(\alpha)}, \partial^-_{p(\beta)}$ and $[p(\alpha):p(\beta)]$ are solid, $\partial^+_{p(a)}, \partial^-_{p(b)}$ and $[p(a):p(b)]$ are dashed. The set $p(J)$ is given by the dotted line.}
    \label{Fig4}
\end{figure}
The following lemma and its proof are visualized in \Cref{Fig4}. The lemma is geometrically intuitive, but the actual mathematical proof convoluted.\\ 
In the proof, we will also use the following observation, which the reader should know before reading the proof. If we have $\alpha<\beta\in[0,2\pi)$ and a supporting hyperplane containing $p(\alpha)$ and $p(\beta)$, then either $p([\alpha,\beta])$ or $p([0,2\pi)\setminus(\alpha,\beta))$ are contained in a face of $K$. 
\begin{Lemma}
Let $I, J\subset [0,2\pi)$ be closed disjoint intervals and $[a,b]\subset [\alpha,\beta]=I\subset[0,2\pi)$ with $a<b$. Then 
  \begin{equation*}
    \textup{dist}(p(J),\partial^+_{p(\alpha)}\cup\partial^-_{p(\beta)}\cup[p(\alpha):p(\beta)])\le\textup{dist}(p(J),\partial^+_{p(a)}\cup\partial^-_{p(b)}\cup[p(a):p(b)]).
  \end{equation*}
\label{LemDis2}
\end{Lemma}

\underline{Proof}:\\
First note that if $p([0,2\pi)\setminus(\alpha,\beta))\subset [p(\alpha):p(\beta)]$, then
  \begin{equation*}
    \textup{dist}(p(J),\partial^+_{p(\alpha)}\cup\partial^-_{p(\beta)}\cup[p(\alpha):p(\beta)])=0
  \end{equation*}
and there is nothing to show. It seems useless to rule out this almost trivial observation, however, we have to get rid of this case later.\\
So let $p([0,2\pi)\setminus(\alpha,\beta))$ not be contained in $[p(\alpha):p(\beta)]$. The following is a general construction for which \Cref{Fig4} is helpful for understanding. The hyperplanes $\partial^+_{p(\alpha)}$ and $\partial^-_{p(\beta)}$ divide $\mathbb{R}^2$ into up to four convex sets, one of which contains $K$ and is given by \Cref{RemLeft} by
  \begin{equation*}
    H=\{z\in\mathbb{R}^2;\ \langle z-p(\beta),ip'_-(\beta)\rangle\ge0,\ \langle z-p(\alpha),ip'_+(\alpha)\rangle\ge0\}.
  \end{equation*}
Furthermore, the hyperplane $[p(\alpha):p(\beta)]$ splits $H$ into two further closed convex sets 
  \begin{equation*}
    \begin{split}
      &H_1=H\cap\{z\in\mathbb{R}^2;\ \langle z-p(\alpha),i(p(\beta)-p(\alpha))\rangle\ge0\},\\
      &H_2=H\cap\{z\in\mathbb{R}^2;\ \langle z-p(\alpha),i(p(\beta)-p(\alpha))\rangle\le0\}.
    \end{split}
  \end{equation*}
Similarly to \Cref{RemLeft}, one sees that $p(J)\subset H_1$ and $p([\alpha,\beta])\subset H_2$. Let us quickly argue why it suffices to show that
  \begin{equation}
    (\partial^+_{p(a)}\cup\partial^-_{p(b)}\cup[p(a):p(b)])\cap \textup{Int}(H_1)=\emptyset.
\label{EqHypDist}
  \end{equation}
Assume that we have shown \Cref{EqHypDist} and that the lemma is false. Then there is a $z\in\partial^+_{p(a)}\cup\partial^-_{p(\beta)}\cup[p(\alpha):p(\beta)]$ and a $t\in J$ such that 
  \begin{equation*}
    \|p(t)-z\|<\textup{dist}(p(J),\partial^+_{p(\alpha)}\cup\partial^-_{p(\beta)}\cup[p(\alpha):p(\beta)]).
  \end{equation*}
However, $\partial H_1\subset \partial^+_{p(\alpha)}\cup\partial^-_{p(\beta)}\cup[p(\alpha):p(\beta)]$, from which it follows that  $z\notin H_1$ and therefore we either have
  \begin{equation*}
    \langle z-p(\beta),ip'_-(\beta)\rangle<0
  \end{equation*}
\textup{or}
  \begin{equation*}
    \langle z-p(\alpha),ip'_+(\alpha)\rangle<0
  \end{equation*}
\textup{or}
  \begin{equation*}
   \langle z-p(\alpha),i(p(\beta)-p(\alpha))<0.
  \end{equation*}
In each of the three cases above there is a $s\in[0,1]$ such that $sz+(1-s)p(t)$ is either  contained in $\partial^+_{p(\alpha)}, \partial^-_{p(\beta)}$ or $[p(\alpha):p(\beta)]$. So we obtain the contradiction
  \begin{equation*}
    \begin{split}
      \|p(t)-z\|&<\textup{dist}(p(J),\partial^+_{p(\alpha)}\cup\partial^-_{p(\beta)}\cup[p(\alpha):p(\beta)])\\
      &\le\|p(t)-(sz+(1-s)p(t))\|=s\|(p(t)-z)\|.
    \end{split}
  \end{equation*}
So it remains to show \Cref{EqHypDist}. Assume there is a $z$ in the left side of \Cref{EqHypDist}. By construction we have that
  \begin{equation*}
    \begin{split}
      &\langle z-p(\alpha), i(p(\beta)-p(\alpha))\rangle>0\ \textup{and}\\
      &\langle z-p(\alpha), ip'_+(\alpha)\rangle>0\ \textup{and}\\
      &\langle z-p(\beta), ip'_-(\beta)\rangle>0,
    \end{split}
  \end{equation*}
since $z\in\textup{Int}(H_1)$, and since $p([a,b])\subset p([\alpha,\beta])\subset H_2$, we also have that
 \begin{equation*}
    \begin{split}
      &\langle p(a)-p(\alpha),i(p(\beta)-p(\alpha))\rangle\le0\ \textup{and}\\
      &\langle p(b)-p(\alpha),i(p(\beta)-p(\alpha))\rangle\le0.
    \end{split}
  \end{equation*}
Thus there is a $w\in\{tz+(1-t)p(a);\ t\in[0,1]\}$ such that 
  \begin{equation*}
    \langle \omega-p(\alpha),i(p(\beta)-p(\alpha))\rangle=0.
  \end{equation*}
We also know that $w\in H\cap[p(\alpha):p(\beta)]$, since $z$ and $p(a)$ are in $H$,.\\
Next we want to show that $\omega\in\{tp(\alpha)+(1-t)p(\beta);\ t\in(0,1)\}$. Let $s\in\mathbb{R}$ such that $w=sp(\alpha)+(1-s)p(\beta)$ and note that
  \begin{equation*}
    0\le\langle w-p(\alpha),ip'_+(\alpha)\rangle=(1-s)\langle p(\beta)-p(\alpha),ip'_+(\alpha)\rangle
  \end{equation*}
and 
\begin{equation*}
    0\le\langle w-p(\beta),ip'_-(\beta)\rangle=s\langle p(\alpha)-p(\beta),ip'_-(\beta)\rangle.
  \end{equation*}
Thus we only have to show that
  \begin{equation*}
    \begin{split}
      &\langle p(\beta)-p(\alpha),ip'_+(\alpha)\rangle>0\ \textup{and}\\
      &\langle p(\alpha)-p(\beta),ip'_-(\beta)\rangle>0.
    \end{split}
  \end{equation*}
Assume not, then 
  \begin{equation*}
    \begin{split}
      &\langle p(\beta)-p(\alpha),ip'_+(\alpha)\rangle=0\ \textup{or}\\
      &\langle p(\alpha)-p(\beta),ip'_-(\beta)\rangle=0,
    \end{split}
  \end{equation*}
since $p(\beta), p(\alpha)\in H$, and hence $p(\alpha)$ and $p(\beta)$ are contained in either $\partial^+_{p(\alpha)}$ or $\partial^-_{p(\beta)}$. According to the note before the theorem, we either have $p([\alpha,\beta])\subset [p(\alpha):p(\beta)]$ or $p([0,2\pi)\setminus(\alpha,\beta))\subset [p(\alpha):p(\beta)]$. The first case implies that $\partial^+_{p(a)}=\partial^-_{p(b)}=[p(a):p(b)]=[p(\alpha):p(\beta)]$ by \Cref{RemHyp}, contradicting the existence of $z$, and the second case contradicts our first assumption of the proof namely that $p([0,2\pi)\setminus(\alpha,\beta))$ is not contained in $[p(\alpha):p(\beta)]$.\\
Thus there is a $s\in[0,1]$ such that $\omega=sp(\alpha)+(1-s)p(\beta)$. We finish the proof with three case distinctions. \\
If $z\in\partial^+_{p(a)}$, then, because of $z\in\textup{Int}(H_1)$, $\partial^+_{p(a)}\neq[p(\alpha):p(\beta)]$ and $w\in\{tz+(1-t)p(a);\ t\in[0,1]\}\subset\partial^+_{p(a)}$, implying that
  \begin{equation*}
    \langle w-p(a),ip'_+(a)\rangle=0
  \end{equation*}
and either 
  \begin{equation*}
    \langle p(\alpha)-p(a),ip'_+(a)\rangle>0
  \end{equation*}
or 
  \begin{equation*}
    \langle p(\beta)-p(a),ip'_+(a)\rangle>0.
  \end{equation*}
However, the equation
  \begin{equation*}
    0=\langle w-p(a),ip'_+(a)\rangle=s\langle p(\alpha)-p(a),ip'_+(a)\rangle+(1-s)\langle p(\beta)-p(a),ip'_+(a)\rangle
  \end{equation*}
shows that $p(\alpha)$ and $p(\beta)$ are not contained in the same half space bounded by $\partial^+_{p(a)}$, contradicting the fact that $\partial^+_{p(a)}$ is a supporting hyperplane.\\
The case $z\in\partial^-_{p(b)}$ is treated analogously by considering $\tilde \omega\in\{tz+(1-t)p(b);\ t\in[0,1]\}$ with
  \begin{equation*}
    \langle \tilde\omega-p(\alpha),i(p(\alpha)-p(\beta))\rangle=0
  \end{equation*}
and doing the above with $\tilde\omega$ instead of $\omega$.\\
If $z\in[p(a):p(b)]$, then $w\in[p(a):p(b)]$. Similarly to \Cref{RemLeft}, we see that
 \begin{equation*}
    \begin{split}
      &\langle p(\alpha)-p(a),i(p(b)-p(a))\rangle\ge0,\\
      &\langle p(\beta)-p(a),i(p(b)-p(a))\rangle\ge0.
    \end{split}
  \end{equation*}
Thus,
  \begin{equation*}
    \begin{split}
      0&=\langle\omega-p(a),i(p(b)-p(a))\rangle\\
      &=s\langle p(\alpha)-p(a),i(p(b)-p(a))\rangle +(1-s)\langle p(\beta)-p(a),i(p(b)-p(a))\rangle
    \end{split}
  \end{equation*}  
implies that $p(\alpha), p(\beta)\in[p(a):p(b)]$ and hence $[p(a):p(b)]=[p(\alpha):p(\beta)]$, contradicting $z\in\textup{Int}(H_1)$. $\hfill\square$

\subsection{Main Theorem}
\label{Section3.3}
Let us quickly recall the setting and notation made at the beginning of this section. $K$ is a compact convex subset of $\mathbb{R}^2$ such that $0\in\textup{Int}(K)$ and $p$ its polar parametrization. For two disjoint closed intervals $[a,b]=I, J\subset[0,2\pi)$ with $a<b$, $g_I$ was the unique composition of a rotation and translation such that $\textup{Im}(g_I(p(a)))=\textup{Im}(g_I(p(b)))=0$ and $\textup{Re}(g_I(p(a)))<\textup{Re}(g_I(p(b)))=0$. In addition, we defined the affine function
  \begin{equation*}
    \tilde g_I=\frac{\textup{Im}(g_I)+\|\textup{Im}(g_I)\|_{p(I),\infty}}{\inf_{p(J)}|\textup{Im}(g_I)|},
  \end{equation*}
whenever $\inf_{p(J)}|\textup{Im}(g_I)|>0$.\\
If $\varphi$ is a u.c.p. map on $C(K)$, then we defined $\varphi(\chi_F)=\mu(K\cap F)$, where $\mu$ is the positive operator valued measure corresponding to $\varphi$.\\
\\
\begin{Lemma}
Let $K$ be as above, $\pi:C(\textup{ex}(K))\to\mathcal{B}(H)$ a unital *-homomorphism and $\phi:C(\textup{ex}(K))\to\mathcal{B}(H)$ a u.c.p. map. If $x\in \textup{ex}(K)$ and $F\subset \textup{ex}(K)\setminus\{x\}$ Borel, then
  \begin{equation*}
    \pi(\chi_{\{x\}})\phi(\chi_F)\pi(\chi_{\{x\}})=0.
  \end{equation*}
\label{LemSingelton}
\end{Lemma}

\underline{Proof}:\\
Let $x\in\textup{ex}(K)$, $E=\textup{ex}(K)\setminus{\{x\}}$ and $\tilde H=\pi(\chi_{\{x\}})(H)$. Assume that $\pi(\chi_{\{x\}})\neq0$. Then the map $P_{\tilde H}\pi|_{\tilde H}$ is a direct sum of copies of the point evaluation $e_x:C(\textup{ex}(K))\to \mathbb{C}$ and hence has the unique extension property by \cite[Proposition 4.4]{WA}. Thus 
  \begin{equation*}
    P_{\tilde H}\phi|_{\tilde H}=P_{\tilde H}\pi|_{\tilde H}\ \textup{on}\ C(\textup{ex}(K)).
  \end{equation*}
If $\mu$ is the positive operator valued measure corresponding to $\phi$, then $P_{\tilde H}\mu|_{\tilde H}$ is the positive operator valued measure corresponding to $P_{\tilde H}\phi|_{\tilde H}$ and by the above equation equal to $\delta_xid_{\tilde H}$. Hence 
  \begin{equation*}
    \pi(\chi_{\{x\}})\phi(\chi_{\{x\}})\pi(\chi_{\{x\}})=P_{\tilde H}
  \end{equation*}
and 
  \begin{equation*}
    \pi(\chi_{\{x\}})\phi(\chi_{E})\pi(\chi_{\{x\}})=P_{\tilde H}-\pi(\chi_{\{x\}})\phi(\chi_{\{x\}})\pi(\chi_{\{x\}})=0.
  \end{equation*}
Now the equality for arbitrary $F\subset E$ Borel follows from 
  \begin{equation*}
    0\le\pi(\chi_{\{x\}})\phi(\chi_{F})\pi(\chi_{\{x\}})\le\pi(\chi_{\{x\}})\phi(\chi_{E})\pi(\chi_{\{x\}})=0.
  \end{equation*}
$\hfill\square$

\begin{Theorem}
Let $K$ be as above, $\pi:C(\textup{ex}(K))\to\mathcal{B}(H)$ a unital *-homomorphism and $\phi:C(\textup{ex}(K))\to\mathcal{B}(H)$ a u.c.p. map such that $\pi=\phi$ on $A(K)$. Let $I, J\subset[0,2\pi)$ be disjoint closed intervals. Then
  \begin{equation*}
    \pi(\chi_{p(I)})\phi(\chi_{p(J)})\pi(\chi_{p(I)})=0.
  \end{equation*}
\label{ThmMain0}
\end{Theorem}

\underline{Proof}:\\
The case that $I$ is a singleton has already been dealt with in \Cref{LemSingelton}. So we can assume that $I=[a,b]$ with $a<b$ and let $\epsilon>0$. By \Cref{LemUpp} there is a partition $a=t_0<t_1<\dots<t_{m}=b$ such that 
  \begin{equation}
    \|\textup{Im}(g_{[t_n,t_{n+1}]})\|_{[t_n,t_{n+1}], \infty}\le\epsilon|p(t_n)-p(t_{n+1})|
\label{EqGeoBou}
  \end{equation}
for $n=0,\dots, m-1$. Denote $[t_n,t_{n+1}]=I_n$. The discussion at the beginning of \Cref{Section3} showed that
  \begin{equation*}
    \begin{split}
          \|\pi(\chi_{p(I)})\phi(\chi_{p(J)})\pi(\chi_{p(I)})\|\le\|&\pi(\chi_{p([t_0,t_{1}])})\phi(\chi_{p(J)})\pi(\chi_{p([t_0,t_{1}])})\|+\\
      &\sum_{n=1}^{m-1}\|\pi(\chi_{p((t_n,t_{n+1}])})\phi(\chi_{p(J)})\pi(\chi_{p((t_n,t_{n+1}])})\|,
    \end{split}
  \end{equation*}
which in turn together with
  \begin{equation*}
    \|\pi(\chi_{p((t_n,t_{n+1}])})\phi(\chi_{p(J)})\pi(\chi_{p((t_n,t_{n+1}])})\|\le\|\pi(\chi_{p(I_n)})\phi(\chi_{p(J)})\pi(\chi_{p(I_n)})\|
  \end{equation*}
implies that
  \begin{equation*}
    \|\pi(\chi_{p(I)})\phi(\chi_{p(J)})\pi(\chi_{p(I)})\|\le\sum_{n=0}^{m-1}\|\pi(\chi_{p(I_n)})\phi(\chi_{p(J)})\pi(\chi_{p(I_n)})\|.
  \end{equation*}
In addition, the discussion showed that
  \begin{equation*}
    \pi(\chi_{p(I_n)})\phi(\chi_{p(J)})\pi(\chi_{p(I_n)})\le\frac{\|\textup{Im}(g_{I_n})\|_{I_n, \infty}}{\inf_{p(J)}|\textup{Im}(g_{I_n})|}id_H
  \end{equation*}
as long as $\inf_{p(J)}|\textup{Im}(g_{I_n})|>0$. The next step is to distinguish the two cases $c>0$ and $c=0$, with
  \begin{equation*}
c=\textup{dist}(p(J),\partial^+_{p(a)}\cup\partial^-_{p(b)}\cup[p(a):p(b)]).
  \end{equation*}
\underline{Step 1:}\\
Assume that $c>0$. By \Cref{LemDis1} and \Cref{LemDis2}, we know that $\textup{inf}_J|\textup{Im}(g_{I_n})|\ge c$ for $n=0,\dots, m-1$. Thus, the above inequalities together with \Cref{EqGeoBou} result in
  \begin{equation}
    \|\pi(\chi_{p(I)})\phi(\chi_{p(J)})\pi(\chi_{p(I)})\|\le \sum_{n=0}^{m-1}\frac{\epsilon}{c}|p(t_n)-p(t_{n+1})|.
\label{EqI}
  \end{equation}
However, the boundary of a convex compact set has finite length \cite[Problem 1.5.1]{VT} and thus, if $L$ is the length of $\partial K$, we conclude that
  \begin{equation*}
    \|\pi(\chi_{p(I)})\phi(\chi_{p(J)})\pi(\chi_{p(I)})\|\le\epsilon\frac{L}{c}.
  \end{equation*}
Since $c$ and $L$ are independent of $\epsilon$, we obtain the theorem.\\
\\
\underline{Step 2}:\\
Assume that $c=0$. Let 
  \begin{equation*}
    \tilde I=\textup{conv}\{p^{-1}(p(I)\cap\textup{ex}(K))\cap[0,2\pi)\}\subset I.
  \end{equation*}
The definition of $\tilde I$ yields that $\textup{ex}(K)\cap p(I\setminus\tilde I)=\emptyset$. Thus the theorem is clear if $\tilde I=\emptyset$ or if $\tilde I$ is a singleton by \Cref{LemSingelton}. So let $\tilde I=[\tilde a,\tilde b]$ with $\tilde a<\tilde b\in[0,2\pi)$. Note that $p(\tilde a), p(\tilde b)\in\textup{ex}(K)$ by definition of $\tilde I$. Let $N\in\mathbb{N}$ such that $\tilde a +1/N<\tilde b -1/N$ and define for $n\ge N$ the intervals $\tilde I_n=[\tilde a+1/n,\tilde b-1/n]$. We claim that
  \begin{equation*}
    \textup{dist}(p(J),\partial^+_{p(\tilde a+1/n)}\cup\partial^-_{p(\tilde b-1/n)}\cup[p(\tilde a+1/n):p(\tilde b-1/n)])>0
  \end{equation*}
for every $n\ge N$. If not, there would be a point $t\in J$ such that either $p(\tilde a+1/n)$ and $p(t)$ or $p(\tilde b-1/n)$ and $p(t)$ are contained in a face of $K$. But this leads to the contradiction that $p(\tilde a)$ or $p(\tilde b)\notin\textup{ex}(K)$. \\
Hence we can use Step 1 and obtain that 
  \begin{equation*}
    \pi(\chi_{p(\tilde I_n)})\phi(\chi_{p(J)})\pi(\chi_{p(\tilde I_n)})=0
  \end{equation*}
for all $n\ge N$. But $\pi(\chi_{p(\tilde I_n)})\to_{SOT}\pi(\chi_{p((\tilde a,\tilde b))})$ and so it also holds that
  \begin{equation*}
     \pi(\chi_{p((\tilde a, \tilde b))})\phi(\chi_{p(J)})\pi(\chi_{p((\tilde a,\tilde b))})=0.
  \end{equation*}
Now $I$ splits up into 
  \begin{equation*}
    I=\{\tilde a\}\cup\{\tilde b\}\cup(\tilde a,\tilde b)\cup (I\setminus \tilde I).
  \end{equation*}
Note that $\pi(\chi_{p(I\setminus \tilde I)})=0$, since $\textup{ex}(K)\cap (I\setminus\tilde I)=\emptyset$, and 
  \begin{equation*}
    \pi(\chi_{\{p(\tilde a)\}})\phi(\chi_{\{p(J)\}})\pi(\chi_{\{p(\tilde a)\}})=0
  \end{equation*}
respectively
  \begin{equation*}
    \pi(\chi_{\{p(\tilde b)\}})\phi(\chi_{\{p(J)\}})\pi(\chi_{\{p(\tilde b)\}})=0 
  \end{equation*}
due to \Cref{LemSingelton}. By using \Cref{LemBrown} once again, we conclude that
  \begin{equation*}
    \pi(\chi_{p(I)})\phi(\chi_{p(J)})\pi(\chi_{p(I)})=0.
  \end{equation*}
$\hfill\square$\\
\ \\
\begin{Corollary}
Under the assumptions of the preceding theorem,
  \begin{equation*}
    \pi(\chi_{p(J)})=\phi(\chi_{p(J)}).
  \end{equation*} 
\label{CorEq}
\end{Corollary}

\underline{Proof}:\\ 
First note that it suffices to show that
  \begin{equation*}
    \phi(\chi_{p(J)})\pi(\chi_{p([0,2\pi)\setminus J)})=0
  \end{equation*}
and 
 \begin{equation*}
    \phi(\chi_{p([0,2\pi)\setminus J)})\pi(\chi_{p(J)})=0
  \end{equation*}
since 
 \begin{equation*}
    \pi(\chi_{p(J)})=\phi(\chi_{p(J)})+\phi(\chi_{p([0,2\pi)\setminus J)})\pi(\chi_{p(J)})-\phi(\chi_{p(J)})\pi(\chi_{p([0,2\pi)\setminus J)}).
  \end{equation*}
Using the Schwarz inequality for u.c.p. maps, we see that
  \begin{equation*}
    \begin{split}
      0&\le(\phi(\chi_{p(J)})\pi(\chi_{p([0,2\pi)\setminus J)}))^*\phi(\chi_{p(J)})\pi(\chi_{p([0,2\pi)\setminus J)})\\
      &\le \pi(\chi_{p([0,2\pi)\setminus J)})\phi(\chi_{p(J)})\pi(\chi_{p([0,2\pi)\setminus J)}),\\
      0&\le(\pi(\chi_{p([0,2\pi)\setminus J)})\pi(\chi_{p(J)}))^*\phi(\chi_{p([0,2\pi)\setminus J)})\pi(\chi_{p(J)})\\
      &\le \pi(\chi_{p(J)})\phi(\chi_{p([0,2\pi)\setminus J)})\pi(\chi_{p( J)}).
    \end{split}
  \end{equation*}
Therefore we only have to show that the right hand sides of the above inequalities are $0$. Let $I_n\subset [0,2\pi)\setminus J$ be an increasing sequence of sets such that every $I_n$ is the union of at most two closed intervals and $\bigcup_n I_n=[0,2\pi)\setminus J$. Then 
  \begin{equation*}
    \begin{split}
      &\pi(\chi_{p(I_n)})\to_{SOT}\pi(\chi_{p([0,2\pi)\setminus J)}),\\
      &\phi(\chi_{p(I_n)})\to_{WOT}\phi(\chi_{p([0,2\pi)\setminus J)})
    \end{split}
  \end{equation*}
and by \Cref{ThmMain0} and \Cref{LemBrown} it also holds that
  \begin{equation*}
    \begin{split}
      &\pi(\chi_{p(J)})\phi(\chi_{p(I_n)})\pi(\chi_{p(J)})=0,\\
      &\pi(\chi_{p(I_n)})\phi(\chi_{p(J})\pi(\chi_{p(I_n)})=0.
    \end{split}
  \end{equation*}
Combing these observations results in
 \begin{equation*}
    \begin{split}
      &\pi(\chi_{p(J)})\phi(\chi_{p([0,2\pi)\setminus J})\pi(\chi_{p(J)})=0,\\
      &\pi(\chi_{p([0,2\pi)\setminus J)})\phi(\chi_{p(J)})\pi(\chi_{p([0,2\pi)\setminus J)})=0.
    \end{split}
  \end{equation*}
$\hfill\square$.

\begin{Theorem}
Let $K\subset \mathbb{R}^2$ be a convex compact set. Then the operator system $A(K)$ is hyperrigid in $C(\textup{ex}(K))$.
\label{ThmMain}
\end{Theorem}

\underline{Proof}:\\
In the discussion at the beginning of \Cref{Section3}, we already saw that the theorem is true if $\textup{Int}(K)=\emptyset$ and that we can reduce the theorem to the case that $0\in\textup{Int}(K)$. \\
Let $\pi$ be a unital *-homomorphism and $\phi$ be a u.c.p. map that agrees with $\pi$ on $A(K)$. By \Cref{CorEq} the set
  \begin{equation*}
    \mathcal{E}=\{E\subset[0,2\pi);\ E\ \textup{Borel}\ , \pi(\chi_{p(E)})=\phi(\chi_{p(E)})\}
  \end{equation*}
contains all closed intervals. It is obvious that the set is closed under complements and also countable unions of disjoint sets. Therefore $\mathcal{E}$ is a Dynkin system. By the Dynkin's $\pi$-$\lambda$ Theorem, we obtain that $\pi(\chi_{E})=\phi(\chi_E)$ for every Borel set $E\subset \textup{ex}(K)$ and hence $\pi=\phi$. $\hfill\square$\\
\\
To conclude the paper, we present a notable application of the main theorem, which yields a generalization of the well known fact that the weak- and strong operator topology agree on the unitary operators.

\begin{Corollary}
For every compact convex $K\subset\mathbb{R}^2$, the weak and strong operator topologies coincide on the set
  \begin{equation*}
    \{T\in\mathcal{B}(H);\ T\ \textup{normal}\ \textup{and}\ \sigma(T)\subset\textup{ex}(K)\}.
  \end{equation*}
\label{CorMain}
\end{Corollary}

\underline{Proof}:\\
Let $T_n, T\in\mathcal{B}(H)$ be normal operators with $\sigma(T_n), \sigma(T)\subset\textup{ex}(K)$ and $T_n\to_{WOT}T$. Define unital *-homomorphims $\pi, \pi_n:C(\textup{ex}(K))\to\mathcal{B}(H)$ with $\pi_n(z)=T_n$ and $\pi(z)=T$. Then, $\pi_n|_{A(K)}\to_{WOT}\pi|_{A(K)}$ and by \Cref{ThmMain} and \cite[Lemma 2.1] {Kleski}, we obtain that $\pi_n\to\pi$ in the strong operator topology.  $\hfill\square$

\bibliography{BibRigid} 
\bibliographystyle{plain}
Fachrichtung Mathematik, Universit\"at des Saarlandes, 66123 Saarbr\"ucken, Germany\\
\textit{Email address:} scherer@math.uni-sb.de

\end{document}